\documentclass{amsart}
\usepackage{amsfonts,amsmath,amssymb}
\usepackage{hyperref}

\newtheorem*{theorem*}{Theorem}
\newtheorem{lemma}{Lemma}[section]
\newtheorem{proposition}[lemma]{Proposition}
\newtheorem{remark}[lemma]{Remark}

\newtheorem{theorem}[lemma]{Theorem}

\newcommand{\cc}{\mathbb{C}}

\newcommand{\Z}{{\mathbb Z}}

\newcommand{\R}{{\mathbb R}}
\newcommand{\C}{{\mathbb C}}

\newcommand{\Fre}{{Fr\'{e}chet \,}}

\newcommand{\cD}{{\mathcal{D}}}

\begin{document}

\author{Avraham Aizenbud} \email{aizenr@yahoo.com}
\address{Faculty of Mathematics and Computer
Science, The Weizmann Institute of Science POB 26, Rehovot 76100,
Israel}
\author{Dmitry Gourevitch}
\email{guredim@yahoo.com}
\address{Faculty of Mathematics and Computer
Science, The Weizmann Institute of Science POB 26, Rehovot 76100,
Israel}
\author{Eitan Sayag}
\email{eitan.sayag@gmail.com}
\address{Institute of Mathematics, The Hebrew University of Jerusalem,
Jerusalem, 91904, Israel}


\title[$(O(V \oplus F), O(V))$ is a Gelfand
pair]{$(O(V \oplus F), O(V))$ is a Gelfand pair\\ for any
quadratic space $V$ over a local field $F$}
\begin{abstract}
Let $V$ be a quadratic space with a form $q$ over an arbitrary
local field $F$ of characteristic different from 2. Let $W=V
\oplus Fe$ with the form $Q$ extending $q$ with $Q(e)=1$. Consider
the standard embedding $\mathrm{O}(V) \hookrightarrow
\mathrm{O}(W)$ and the two-sided action of $\mathrm{O}(V)\times
\mathrm{O}(V)$ on $\mathrm{O}(W)$.

In this note we show that any $\mathrm{O}(V) \times
\mathrm{O}(V)$-invariant distribution on $\mathrm{O}(W)$ is
invariant with respect to transposition. This result was earlier
proven in a bit different form in \cite{vD} for $F=\R$, in
\cite{Apa-vD} for $F=\C$ and in \cite{BvD} for $p$-adic fields.
Here we give a different proof.

Using results from \cite{AGS}, we show that this result on
invariant distributions implies that the pair $(O(V),O(W))$ is a
Gelfand pair. In the archimedean setting this means that for any
irreducible admissible smooth \Fre representation $(\pi,E)$ of
$\mathrm{O}(W)$ we have $ dim Hom_{\mathrm{O}(V)}(E,\cc) \leq 1.$

A stronger result for $p$-adic fields is obtained in \cite{AGRS}.
\end{abstract}
\maketitle \tableofcontents
\section{Introduction}
Let $F$ be a local field of characteristic different from 2.

Let $(W,Q)$ be a quadratic space defined over $F$ and fix $e \in
W$ a unit vector. Consider the quadratic space $V={e}^{\bot}$ with
$q=Q|_{V}$. Define the standard imbedding $\mathrm{O}(V)
\hookrightarrow \mathrm{O}(W)$ and consider the two-sided action
of $\mathrm{O}(V) \times \mathrm{O}(V)$ on $\mathrm{O}(W)$ defined
by
$(g_1,g_2)h:=g_1 hg_2^{-1}$. We also consider the anti-involution
$\tau$ of $O_{Q}$ given by $\tau(g)=g^{-1}$. In this paper we
prove the following theorem
\begin{theorem*} [A]
Any $\mathrm{O}(V) \times \mathrm{O}(V)$ invariant distribution on
$\mathrm{O}(W)$ is invariant under $\tau$.
\end{theorem*}
This theorem has the following corollary in representation theory.

\begin{theorem*}[B]
Let $(\pi,E)$ be an irreducible admissible representation of
$\mathrm{O}(W)$. Then
$$dim Hom_{\mathrm{O}(V)}(E,\cc) \leq 1$$
\end{theorem*}

Here admissible representation refers to the usual notion in the
non-archimedean case and to the notion of admissible smooth \Fre
representation in the archimedean setting.

Our proof for the archimedean and non-archimedean case is uniform,
except at one point where the archimedean case requires an extra
analysis of a certain normal bundle (see lemma \ref{non uniform}).


\begin{remark}
We note that a related result for unitary representations of
$SO(V,Q)$ is proved in \cite{BvD} (for $p$-adic fields) and
in \cite{vD} (for the real numbers). In fact, the proof given in
those papers implies also theorem \textbf{A}. Also, an analogous
theorem for unitary groups is proven in \cite{vD2}.
\end{remark}

\subsection*{Acknowledgements}
We thank Prof. Gerrit van Dijk for pointing out to us that the
arguments of \cite{vD}, \cite{Apa-vD} and \cite{BvD} give a proof
of theorem \textbf{A}. We also thank Dr. Sun Binyong for finding a
mistake in the previous version of this note. Finally, we would
like to thank the referee for useful remarks.

\section{From Invariant distributions to Representation theory}
In this section we recall a technique due to Gelfand and Kazhdan
which allows to deduce theorem {\bf B} from theorem {\bf A}.


 Recall the following theorem (\cite{AGS})

\begin{theorem}\label{GKThomasCW}
Let $H \subset G$ be reductive groups and let $\tau$ be an
involutive anti-automorphism of $G$ and assume that $\tau(H)=H$.
Suppose $\tau(T)=T$ for all bi $H$-invariant distributions
\footnote{In fact it is enough to check this only for Schwartz
distributions.} on $G$. Then for any irreducible admissible
representation $(\pi,E)$ of $G$ we have

$$dim Hom_{H}(E,\cc) \cdot dim Hom_{H}(\widetilde{E},\cc)\leq
1,$$

where $\widetilde{E}$ denotes the smooth contragredient
representation.

\end{theorem}

Note that in the non-archimedean case the same result is proven in
\cite{Prasad}.

 To finish the deduction of theorem \textbf{B} from theorem \textbf{A} we will show that
\begin{theorem}\label{O}
Let $(\pi,E)$ be an irreducible admissible representation of
$G=O(V)$. Then $\widetilde{E} \cong E$ and in particular
$$dim Hom_{H}(E,\cc) =dim Hom_{H}(\widetilde{E},\cc)$$
\end{theorem}
For the proof we recall proposition I.2 (chapter 4) from
\cite{MVW}:

\begin{proposition}
Let $V$ be a quadratic space and let $g \in O(V)$. Then $g$ is
conjugate to $g^{-1}$.
\end{proposition}

\begin{proof}[Proof of Theorem \ref{O}]
For non-archimedean fields this is a theorem from \cite{MVW} page
91.  For archimedean fields we use the Harish-Chandra regularity
theorem and the proposition that any element in $g \in O(V)$ is
conjugate in $O(V)$ to $g^{-1}$. Thus, the characters of $E$ and
$\widetilde{E}$ are the same and hence $\widetilde{E} \cong E$.
\end{proof}

\begin{remark}
A related result for the groups $SO(V)$ can be found in \cite{GP},
proposition 5.3.
\end{remark}

\section{Basic Results on Invariant distributions}  \label{ProofLoc}

In this paper we consider distributions over $l$-spaces and over
smooth manifolds. $l$-spaces are locally compact totally
disconnected topological spaces (see \cite{BZ}, section 1).

For $X$ a smooth manifold or an $l$-space we denote by $\cD(X)$
the space of distributions on $X$. When $X$ is an $l$-space this
means that $\cD(X)=S(X)^{*}$ where $S(X)$ is the space of locally
constant functions with compact support on $X$. For smooth $X$, we
let $\cD(X)=C_{c}^{\infty}(X)^{*}$.

The basic tools to study invariant distributions on a $G$-space
$X$ are Bruhat filtration, Frobenuis reciprocity (\cite{BZ},
\cite{Bar} and \cite{AGS}) and the Bernstein's localization
principle (\cite{Ber} and \cite{AG}). Let us remind the
statements.

For the simplicity of formulation we provide, for each principle,
two versions: for $l$-spaces and for smooth manifolds.

\subsection{Bruhat Filtration}
Although we will not need the non-archimedean version of this
principle, we formulate it for completeness. It is a simple
consequence of proposition 1.8 in \cite{BZ}.

\begin{theorem}\label{Filt_nonarch}
Let an $l$-group $G$ act on an $l$-space $X$. Let $X =
\bigcup_{i=0}^l X_i$ be a $G$-invariant stratification of $X$. Let
$\chi$ be a character of $G$. Suppose that $\cD(X_i)^{G,\chi}=0$.
Then $\cD(X)^{G,\chi}=0$.
\end{theorem}

To formulate the archimedean version we let $X$ be a smooth
manifold and $Y \subset X$ a smooth submanifold. We remind the
definition of the conormal bundle $CN_{Y}^{X}$. For this denote by
$T_{X}$ the tangent bundle of $X$ and by
$N_Y^X:=(T_X|_Y)/T_Y $ the normal bundle to $Y$ in $X$. The
conormal bundle is defined by $CN_Y^X:=(N_Y^X)^*$. Denote by
$Sym^k(CN_Y^X)$ the k-th symmetric power of the conormal bundle.

\begin{theorem} \label{Filt}
Let a real reductive group $G$ act on a smooth affine real
algebraic variety $X$. Let $X = \bigcup_{i=0}^l X_i$ be a smooth
$G$-invariant stratification of $X$. Let $\chi$ be an algebraic
character of $G$. Suppose that for any $k \in \Z_{\geq 0}$ and any
$0 \leq i \leq l$ we have $\cD(X_i,Sym^k(CN_{X_i}^X))^{G,\chi}=0$.
Then $\cD(X)^{G,\chi}=0$.
\end{theorem}
For proof see \cite{AGS}, section B.2.

\subsection{Frobenius reciprocity}

For $l$-space, the following version of Frobenius reciprocity is
proven in \cite{Ber}:

\begin{theorem}[Frobenius reciprocity] \label{padic-Frob}
Let a unimodular $l$-group $G$ act transitively on an $l$-space
$Z$. Let $\varphi:X \to Z$ be a $G$-equivariant continuous map.
Let $z\in Z$. Suppose that its stabilizer $\mathrm{Stab}_G(z)$ is
unimodular. Let $X_z$ be the fiber of $z$. Let $\chi$ be a
character of $G$. Then $\cD(X)^{G,\chi}$ is canonically isomorphic
to $\cD(X_z)^{\mathrm{Stab}_G(z),\chi}$.
\end{theorem}

An archimedean version is considered in \cite{Bar}. Here is a
slight generalization (see \cite{AGS}):

\begin{theorem}[Frobenius reciprocity] \label{Frob}
Let a unimodular Lie group $G$ act transitively on a smooth
manifold $Z$. Let $\varphi:X \to Z$ be a $G$-equivariant smooth
map. Let $z\in Z$. Suppose that its stabilizer
$\mathrm{Stab}_G(z)$ is unimodular. Let $X_z$ be the fiber of $z$.
Let $\chi$ be a character of $G$. Then $\cD(X)^{G,\chi}$ is
canonically isomorphic to $\cD(X_z)^{\mathrm{Stab}_G(z),\chi}$.
Moreover, for any $G$-equivariant bundle $E$ on $X$,
$\cD(X,E)^{G,\chi}$ is canonically isomorphic to
$\cD(X_z,E|_{X_z})^{\mathrm{Stab}_G(z),\chi}$.
\end{theorem}

\subsection{Bernstein's Localization principle}

For $l$-spaces it is taken from \cite{Ber}:

\begin{theorem}[Localization principle] \label{padic-Localization}
Let $X$ and $T$ be $l$-spaces and $\phi:X \to T$ be a continuous
map. Let an $l$-group $G$ act on $X$ preserving the fibers of
$\phi$. Let $\chi$ be a character of $G$. Suppose that for any
$t\in T$, $\cD(\phi^{-1}(t))^{G,\chi}=0$. Then
$\cD(X)^{G,\chi}=0$.
\end{theorem}

For real smooth algebraic varieties, the following theorem is
proven in \cite{AG}, Corollary A.0.3:

\begin{theorem}[Localization principle]
Let a real reductive group $G$ act on a smooth affine real
algebraic variety $X$. Let $Y$ be a smooth real algebraic variety
and $\phi:X \to Y$ be an algebraic $G$-invariant submersion.
Suppose that for any $y \in Y$ we have
$\cD(\phi^{-1}(y))^{G,\chi}=0$. Then $\cD(X)^{G,\chi}=0$.
\end{theorem}


%

\section{Proof of Theorem {\bf A}}

Recall the setting. $(W,Q)$ is a quadratic space over $F$,  $e \in
W$ with $Q(e)=1$. Also $(V,q)$ is defined by $V={e}^{\bot}$ and
$q=Q|_{V}$.

We need some further notations.

\begin{itemize}
\item
$O_{q}=O(V,q)$ is the group of isometries of the quadratic space
$(V,q)$.

\item
 $G_{q}=O(V,q) \times O(V,q)$.
\item
 $\Delta: O_{q} \to
G_{q}$ the diagonal. $H_{q}=\Delta(O_{q}) \subset G_{q}$.

\item
$\tau(g_{1},g_{2})=(g_{2},g_{1})$.

\item
$\widetilde{G_{q}}=G_{q} \rtimes \{1,\tau\}$, same for
$\widetilde{H_{q}}$

\item
$\chi : \widetilde{G_{q}} \to \{+1,-1\}$ the non trivial character
with $\chi(G_{q})=1$.

\item
$\widetilde{G_{Q}}$ acts on $O_{Q}$ by
$(g_{1},g_{2})x=g_{1}xg_{2}^{-1}$ and $\tau(x)=x^{-1}$.

\end{itemize}

Clearly Theorem {\bf A} follows from the following theorem:

\begin{theorem} \label{Reform}

$\cD(O_{Q})^{\widetilde{G_{q}},\chi}=0$

\end{theorem}

\subsection{Proof of theorem \ref{Reform}}
We denote by $\Gamma=\{w \in W: Q(w)=1\}$. Note that by Witt's
theorem $\Gamma$ is an $O_{Q}$ transitive set and therefore
$\Gamma \times \Gamma$ is a transitive $\widetilde{G}_{Q}$ set
where the action of $G_{Q}$ is the standard action on $W \oplus W$
and $\tau$ acts by flip.

Applying Frobenuis reciprocity (\ref{padic-Frob}, \ref{Frob}) to
projections of $O_Q \times \Gamma \times \Gamma$ first on $\Gamma
\times \Gamma$ and then on $O_Q$ we have
$$\cD
(O_Q)^{\widetilde{G_{q}},\chi}=\cD(O_{Q} \times \Gamma \times
\Gamma)^{\widetilde{G_{Q}},\chi}$$ and also that

$$\cD(O_{Q} \times \Gamma
\times \Gamma)^{\widetilde{G_{Q}},\chi}=\cD(\Gamma \times
\Gamma)^{\widetilde{H_{Q}},\chi}$$

In what follows we will abuse notation and write $Q(u,v)$ for the
bilinear form defined by $Q$. Define a map $D: \Gamma \times
\Gamma \to Z$ where $Z=\{(v,u) \in W \oplus W: Q(v,u)=0 ,
Q(v+u)=4\}$ by
$$D(x,y)=(x+y,x-y).$$ $D$ defines an $\widetilde{G}_{Q}$-equivariant
homeomorphism and thus we need to show that
$$\cD(Z)^{\widetilde{H_{Q}},\chi}=0$$
Here, the action of $G_{Q}$ on $Z \subset W \oplus W$ is the
restriction of its action on $W \oplus W$ while the action of
$\tau$ is given by $\tau(v,u)=(v,-u).$

Now we cover $Z = U_{1} \cup U_{2}$ where
$$U_{1}=\{(v,u) \in Z: Q(v) \ne 0\}$$ and $$U_{2}=\{(v,u) \in Z: Q(u) \ne 0\}$$
We will show $\cD(U_{1})^{\widetilde{H_{Q}},\chi}=0$, and the
proof for $U_{2}$ is analogous. This will finish the proof.

\begin{lemma}\label{non uniform}
$\cD(U_{1})^{\widetilde{H_{Q}},\chi}=0$ \end{lemma}

\begin{proof}[Proof for non-archimedean $F$.]

Consider $\ell_{1}:U_{1} \to F-\{0\}$ defined as
$\ell_{1}(v,u)=Q(v)$. By the localization principle, it is enough
to show $\cD(U_{1}^{\alpha})^{\widetilde{H_{Q}},\chi}=0$ where
$U_{1}^{\alpha}=\ell_{1}^{-1}(\alpha)$, for any $\alpha \in
F-\{0\}$. But

$$U_{1}^{\alpha}=\{(v,u)| Q(v)=\alpha, Q(u)=4-\alpha, Q(v,u)=0\}$$

Let $W^{\alpha}=\{w \in W| Q(w)=\alpha\}$ and let
$p_1:U_{1}^{\alpha} \to W^{\alpha}$ be given by $p_1(v,u)=v$.

On $W^{\alpha}$ our group acts transitively. Fix a vector $v_{0}
\in W^{\alpha}$.\\ Denote $H(v_0):= H_{(Q|_{v_{0}^{\bot}})}$ and
$\widetilde{H}(v_0):= \widetilde{H}_{(Q|_{v_{0}^{\bot}})}$.

The stabilizer in $\widetilde{H}_{Q}$ of $v_{0}$ is
$\widetilde{H}(v_0)$. The fiber $p_1^{-1}(v_{0})=\{a \in
{v_{0}^{\bot}}| Q(a)=4-\alpha\}.$

Frobenius reciprocity implies that

$$\cD(U_{1}^{\alpha})^{\widetilde{H_{Q}},\chi}=\cD(p_1^{-1}(v_{0}))^{\widetilde{H}(v_0),\chi}$$
But clearly $\cD(p_1^{-1}(v_{0}))^{\widetilde{H}(v_0),\chi}=0$ as
$-Id \in H(v_0)$.
\end{proof}

\begin{proof}[Proof for archimedean $F$.]
Now let us consider the archimedean case. Define $U:=\{(v,u) \in
U_1 | u \neq 0 \}$. Note that the map $\ell _1|_U$ is a
submersion, so the same argument as in the non-archimedean case
shows that $\cD(U)^{\widetilde{H_{Q}},\chi}=0$. Let $Y:=\{(v\in W
| Q(v) = 4\} \times \{0\}$ be the complement to $U$ in $U_1$.
By theorem \ref{Filt}, it is enough to prove
$\cD(Y,Sym^k(CN_Y^{U_1}))^{\widetilde{H_{Q}},\chi}=0$.

Note that the action of $\widetilde{H}_{Q}$ on $Y$ is transitive,
and fix a point $(v,0) \in Y$. The stabilizer in
$\widetilde{H}_{Q}$ of $(v,0)$ is $\widetilde{H}(v)$, and the
normal space to $Y$ at $(v,0)$ is $v^{\bot}$. So Frobenius
reciprocity (theorem \ref{Frob}) implies that

$$\cD(Y,Sym^k(CN_Y^{U_1}))^{\widetilde{H_{Q}},\chi}=Sym^k(v^{\bot})^{\widetilde{H}(v),\chi}$$
But clearly $Sym^k(v^{\bot})^{\widetilde{H}(v),\chi}=0$ as $-Id
\in H(v)$.
\end{proof}


\begin{thebibliography}{99}
\bibitem[AvD]{Apa-vD} Sofia Aparicio and G. van Dijk: {\it Complex generalized Gelfand pairs}, Tambov University Reports
(2006).

\bibitem[\href{http://arxiv.org/abs/0803.3395}{AG}]{AG} A. Aizenbud,  D. Gourevitch: {\it
Generalized Harish-Chandra Descent and applications to Gelfand
Pairs}. arXiv:0803.3395v7  [math.RT].

\bibitem[\href{http://arxiv.org/pdf/0709.1273v4}{AGS}]{AGS} A. Aizenbud,  D. Gourevitch, E. Sayag : {\it
$(GL_{n+1}(F),GL_n(F))$ is a Gelfand pair for any local field
$F$}, postprint: arXiv:0709.1273v4[math.RT]. Originally published
in: Compositio Mathematica, \textbf{144} , pp 1504-1524 (2008),
doi:10.1112/S0010437X08003746.

\bibitem[\href{http://arxiv.org/pdf/0709.4215v1}{AGRS07}]{AGRS} A. Aizenbud,  D. Gourevitch, S. Rallis, G. Schiffmann, {\it
Multiplicity One Theorems}, arXiv:0709.4215v1 [math.RT], To appear
in the Annals of Mathematics.

\bibitem[\href{http://projecteuclid.org/DPubS?service=UI&version=1.0&verb=Display&handle=euclid.annm/1061030449}{Bar}]{Bar} E.M. Baruch:
{\it A proof of Kirillov's conjecture.} Annals of Mathematics,
\textbf{158}, 207–252 (2003).

\bibitem[\href{http://www.math.tau.ac.il/~bernstei/Publication_list/publication_texts/Bernstein-P-invar-SLN.pdf}{Ber}]{Ber} J. Bernstein:
{\it $P$-invariant Distributions on $\mathrm{GL}(N)$ and the
classification of unitary representations of $\mathrm{GL}(N)$
(non-archimedean case)} Lie group representations, II (College
Park, Md., 1982/1983), 50--102, Lecture Notes in Math.,
\textbf{1041}, Springer, Berlin (1984).

\bibitem[BvD]{BvD} E. E H. Bosman and G. Van Dijk: {\it A New Class of Gelfand Pairs.}
Geometriae Dedicata 50: 261-282, 1994. 261 @ 1994 KluwerAcademic
Publishers. Printed in the Netherlands.

\bibitem[\href{http://www.math.tau.ac.il/~bernstei/Publication_list/publication_texts/B-Zel-RepsGL-Usp.pdf}{BZ}]{BZ} J. Bernstein, A.V.
Zelevinsky: {\it Representations of the group $\mathrm{GL}(n, F)$,
where F is a local non-Archimedean field,} Uspekhi Mat. Nauk
\textbf{10}, No.3, 5-70 (1976).

\bibitem [GP]{GP} Gross, Benedict H. and Prasad, Dipendar: \textit{On the decomposition of a representation of $SO_n$ when restricted to $SO_{n-1}$}, Can. J. Math, \textbf{44} (1992), no 5, pp 974-1002

 \bibitem[MVW]{MVW} M\oe glin, Colette and  Vigneras, Marie-France and  Waldspurger, Jean-Loup: \textit{Correspondances de Howe sur un corps $p$-adique. (French)
 [Howe correspondences over a $p$-adic field]} Lecture Notes in Mathematics, \textbf{1291},
 Springer-Verlag, Berlin, 1987. viii+163 pp. ISBN: 3-540-18699-9

\bibitem[Pra]{Prasad} D. Prasad: {\it Trilinear forms
for representations of $GL_{2}$ and local $\epsilon$ factors}.
Compositio Mathematica Tome 75, N.1 ,1990, page 1-46
\bibitem[Ser]{Ser} J.P. Serre: {\it Lie Algebras and Lie Groups} Lecture Notes in Mathematics \textbf{1500}, Springer-Verlag, New York, (1964).


\bibitem[vD]{vD} G. van Dijk: {\it On a class of generalized Gelfand
pairs} Math. Z. 193, 581-593 (1986).

\bibitem[vD2]{vD2} G. van Dijk: {\it $(U(p, q), U(p - 1, q))$ is a generalized Gelfand
pair}. Preprint.

\end{thebibliography}
\end{document}